\documentclass{article}
\usepackage{graphicx} 
\usepackage{hyperref}
\usepackage{amssymb}
\usepackage{amsmath}
\usepackage{amsthm}
\newtheorem{theorem}{Theorem}

\title{The Maximum of the Volume of a Cevian Simplex and its Parts}
\author{Yagub Aliyev}
\date{April 2026}

\begin{document}

\maketitle

\begin{abstract}
The cevian triangle corresponding to an interior point $M$ of a triangle is the triangle determined by
the feet of the three cevians concurrent at $M$. It is known that the area of the cevian triangle for an interior point $M$ of a triangle is at most $\frac{1}{4}$ of the area of the triangle, with maximum attained when $M$ is the
triangle's centroid. This can be generalized from triangles to $n$-dimensional simplices,
with $\frac{1}{4}$ replaced by $\frac{1}{n^n}$, using barycentric coordinates. We also use this method to solve two optimization problems about the parts of this simplex.
\end{abstract}


\medskip


\medskip

The area-ratio of any cevian triangle and
the base triangle is known to be at most $\frac{1}{4}$ (see \cite{bog}), with this maximum attained when the cevian
triangle is the median triangle. This also follows from Möbius' theorem (see \cite{balk}, p. 95-96, \cite{ali}, \cite{moeb}, p. 198). This asserts that if $A_1N_1$, $A_2N_2$, and $A_3N_3$, are the cevians of a triangle concurrent at the point $M$, and $p$, $q$, $r$, and $x$ are the areas of the triangles $A_1N_2N_3$, $A_2N_3N_1$, $A_3N_1N_2$, and $N_1N_2N_3$, respectively, then $4pqr=x^2(p+q+r+x)$. Indeed, if the area of triangle $A_1A_2A_3$ is denoted by $S$, then by the Arithmetic Mean - Geometric Mean (AM-GM) inequality
$$
\frac{S}{4}=\frac{p+q+r+x}{4}\ge\sqrt[4]{pqrx}=\sqrt[4]{\frac{1}{4}x^3S}.
$$
By simplifying the inequality of the first and last
terms above, we obtain $x\le\frac{1}{4}S$. Equality holds if and only if $p = q = r = x$;
that is, if and only if the cevian triangle is the same as the median
triangle. The next theorem generalizes this result to higher dimensions \cite{reut}, \cite{klam}, \cite{mit}, p. 502. Throughout this paper, $[X_1X_2\ldots X_{n+1}]$ denotes the volume of an $n$-simplex $X_1X_2\ldots X_{n+1}$.
\begin{theorem} \textnormal{(\cite{reut}, \cite{klam}, \cite{mit}, p. 502)} 
Let $N_1N_2\ldots N_{n+1}$ be the cevian simplex
corresponding to an interior point $M$ of an $n$-simplex $A_1A_2\ldots A_{n+1}$.
Then
$$[N_1N_2\ldots N_{n+1}]\le\frac{[A_1A_2\ldots A_{n+1}]}{n^n},\eqno(1)$$
with equality if and only if $M$ is the centroid
of $A_1A_2\ldots A_{n+1}$.
\end{theorem}
\begin{proof}
Let $(\lambda_1,\lambda_2,\ldots,\lambda_{n+1})$ be the barycentric coordinates of $M$ with respect to the simplex $A_1A_2\ldots A_{n+1}$. If $M$ is on $A_1A_2\ldots A_{n+1}$, then $\lambda_i=0$ or $1$ for some $i$ and the inequality is trivially true because the left side of it vanish. So, in the following, we will assume that $M$ is in the interior of $A_1A_2\ldots A_{n+1}$, which means that $0<\lambda_i<1$ for all $i$. Let $R_i$ denote the distance from $M$ to $A_i$, and let $s_i$ denote the length of $MN_i$. It is well known that $\lambda_i=\frac{s_i}{R_i+s_i}$ and $\lambda_1+{\lambda_2}+\cdots+\lambda_{n+1}=1$ (see \cite{nguy1}; \cite{balk}, p. 124-126, \cite{mark}, \cite{samet}). So, $\frac{s_i}{R_i}=\frac{\lambda_i}{1-\lambda_i}$, and therefore

\begin{flalign*}
& [MN_1N_2\ldots N_{n}]=[MA_1A_2\ldots A_{n}]\cdot \prod_{i=1}^n\frac{\lambda_i}{1-\lambda_i} \\
&=[A_1A_2\ldots A_{n+1}]\cdot \lambda_{n+1}\cdot \prod_{i=1}^n\frac{\lambda_i}{1-\lambda_i}\tag{2} \\
& =[A_1A_2\ldots A_{n+1}]\cdot \left(\prod_{i=1}^{n+1}{\lambda_i}\right)\cdot \left(\prod_{i=1}^n\frac{1}{\sum_{j\ne i}\lambda_j}\right) \\
& \le [A_1A_2\ldots A_{n+1}]\cdot \left(\prod_{i=1}^{n+1}{\lambda_i}\right)\cdot \left(\prod_{i=1}^n\frac{1}{n\sqrt[n]{\prod_{j\ne i}\lambda_j}}\right) \\
& =[A_1A_2\ldots A_{n+1}]\cdot \left(\prod_{i=1}^{n+1}{\lambda_i}\right)\cdot \frac{1}{n^n}\cdot \left(\frac{1}{\sqrt[n]{\left(\prod_{i=1}^n\lambda_i^{n-1}\right)\lambda_{n+1}^n}}\right) \\
& =\frac{[A_1A_2\ldots A_{n+1}]}{n^n}\sqrt[n]{\prod_{i=1}^n\lambda_i}\\
&\le\frac{[A_1A_2\ldots A_{n+1}]}{n^n}\cdot \frac{\sum_{i=1}^n\lambda_i}{n}  =\frac{[A_1A_2\ldots A_{n+1}]}{n^{n+1}}\cdot (1-\lambda_{n+1}),
\end{flalign*}
where both inequalities follow from the AM-GM
inequality. Similarly, for each $i=1,2,\ldots,n+1$,
$$[MN_1N_2\ldots N_{i-1}N_{i+1}\ldots,N_{n+1}]\le\frac{[A_1A_2\ldots A_{n+1}]}{n^{n+1}}\cdot (1-\lambda_{i}).$$
Therefore
$$[N_1N_2\ldots N_{n+1}]=\sum_{i=1}^{n+1}[MN_1N_2\ldots N_{i-1}N_{i+1}\ldots,N_{n+1}]=$$

$$\le\frac{[A_1A_2\ldots A_{n+1}]}{n^{n+1}}\cdot \sum_{i=1}^{n+1}(1-\lambda_{i})=\frac{[A_1A_2\ldots A_{n+1}]}{n^{n}},$$
as needed.
\end{proof}
We can now prove the following surprising result using the formula (2).
\begin{theorem} With the notation from Theorem 1,
$$[MN_1N_2\ldots N_{n}]\le\frac{(n-1)^2}{(n-\theta_n)^{n+3}}\cdot[A_1A_2\ldots A_{n+1}],\eqno(3)
$$
where $\theta_n=\frac{n+1-\sqrt{n^{2}+2 n-3}}{2}$. Equality
happens when $M$ has barycentric coordinates $(\theta_n,\theta_n,\ldots,\theta_n,1-n\theta_n)$ in the given order.
\end{theorem}
\begin{proof}
Since the function $F(\lambda_1,\lambda_2,\ldots,\lambda_{n+1})=\lambda_{n+1}\prod_{i=1}^n\frac{\lambda_i}{1-\lambda_i}$ is continuous on the compact simplex
$A_1A_2\ldots A_{n+1}$, $F$ attains its extremum values in this simplex. Note that on the boundary of the simplex
$A_1A_2\ldots A_{n+1}$ the function $F$ needs to be defined as $F = 0$ to guarantee its continuity for these points, too. As the
minimum of $F = 0$ is attained on the faces of the simplex, the maximum occurs in the interior of the simplex. By the Method of Lagrange Multipliers, for a fixed ${\lambda_i}+{\lambda_j}$ the product $\frac{\lambda_i}{1-\lambda_i}\cdot \frac{\lambda_j}{1-\lambda_j}$ reaches the maximum when $\lambda_i={\lambda_j}$. Therefore, for fixed $\sum_{i=1}^n{\lambda_i}=nx$ the product $\prod_{i=1}^n\frac{\lambda_i}{1-\lambda_i}$ reaches the maximum when $\lambda_1={\lambda_2}=\cdots=\lambda_n=x$. Then $\lambda_{n+1}=1-nx$ and $$\lambda_{n+1}\cdot\prod_{i=1}^n\frac{\lambda_i}{1-\lambda_i}=\left(\frac{x}{1-x}\right)^{n}\cdot \left(1-nx\right)=:f(x).$$
We obtain the maximum of $f$ in the interval $\left(0,\frac{1}{n}\right)$ by standard single variable calculus. The derivative $$f'(x)=\left(\frac{x}{1-x}\right)^{n} \cdot\frac{n \left(x^{2}-(n+1)x+1\right)}{x \left(1-x\right)},$$leads
to the critical point $x=\theta_n$ being the only global maximum point of $f$
in the reference interval. Therefore, $$f(x)\le f \left(\theta_n\right)=\frac{(n-1)^2}{(n-\theta_n)^{n+3}},$$
as needed. The last expression is obtained from the earlier definition
of $f$, using three equivalent forms of the quadratic equation satisfied by $\theta_n$: $1-n\theta_n=\theta_n(1-\theta_n)$, $1-\theta_n=\theta_n(n-\theta_n)$, and $(1-\theta_n)^2=(n-1)\theta_n$.
\end{proof}
When $n=2$, Theorem 2 gives the following result for triangles
$$[MN_1N_2]\le
\frac{32}{\left(\sqrt{5}+1\right)^{5}}\cdot[A_1A_2A_{3}],
$$
which is in perfect agreement with Theorem 1 in \cite{ox}. The same article observed that the
above fraction is $\frac{1}{\phi^5}$, where $\phi$ is the golden ratio.
When $n=3$, we obtain the solution in \cite{nguy0} (see also \cite{nguy2}) to the problem about tetrahedra:
$$[MN_1N_2N_3]\le
\frac{4}{\left(1+\sqrt{3}\right)^{6}}\cdot[A_1A_2A_{3}A_{4}].
$$
Since its barycentric coordinates involve only quadratic irrationalities, the point $M$ maximizing the volume of $MN_1N_2\ldots N_{n}$ is constructible with straightedge and compass for
each $n$. Note that $\theta_2=\frac{1}{\phi^2}$, $\theta_3=\tan{15^\circ}$, and in general $\theta_n$ ($n=1,2,\ldots$) can be written as
$$
\theta_n=e^{-\cosh^{-1}{\frac{n+1}{2}}}=\frac{1}{n+1-\frac{1}{n+1-\frac{1}{n+1-\cdots}}},
$$
which shows a remarkable resemblance to the metallic means $\phi_n=\frac{n+\sqrt{n^2+4}}{2}$ ($n=1,2,\ldots$), also known as the silver means (see for example \cite{wal}, p. 23 for the continued fraction representation and \cite{gil}, \cite{spi}, \cite{weis} for the metallic/silver means):
$$
\phi_n=e^{\sinh^{-1}{\frac{n}{2}}}=n+\frac{1}{n+\frac{1}{n+\frac{1}{n+\cdots}}}.
$$
We prove the following result using the same method.
\begin{theorem}
With the notation from Theorem 1,
$$\frac{[N_1N_2\ldots N_{n}N_{n+1}]-[MN_1N_2\ldots N_{n}]}{[A_1A_2\ldots A_{n+1}]}
\le\frac{\left(n-1\right)^{n-1}\left(1-n\omega_n\right)}{\left(n^{2}-n+1-n\omega_n\right)^{n-1}},$$
where $\omega_n=\frac{n^{2}-\sqrt{n^{4}-4 n^{2}+4 n}}{2 n}$. Equality
happens when $M$ has barycentric coordinates $(\omega_n,\omega_n,\ldots,\omega_n,1-n\omega_n)$ in the given order.
\end{theorem}

\subsection*{Acknowledgement.} The paper was financially supported by the ADA University Faculty Research and Development Fund.

\end{document}